\documentclass[12pt]{article}
\usepackage{graphicx}
\usepackage{amsmath,amsthm,amssymb,enumerate}
\usepackage{euscript,mathrsfs}
\usepackage[left=2cm,right=2cm,top=3.5cm,bottom=3.5cm]{geometry}
\usepackage{color}
\catcode`\@=11 \@addtoreset{equation}{section}

\catcode`\@=12

\allowdisplaybreaks

\newtheorem{Theorem}{Theorem}[section]
\newtheorem{Proposition}[Theorem]{Proposition}
\newtheorem{Lemma}[Theorem]{Lemma}
\newtheorem{Corollary}[Theorem]{Corollary}

\theoremstyle{definition}
\newtheorem{Definition}[Theorem]{Definition}

\newtheorem{Remark}[Theorem]{Remark}

\newcommand{\bTheorem}[1]{
\begin{Theorem} \label{T#1} }
\newcommand{\eT}{\end{Theorem}}

\newcommand{\bProposition}[1]{
\begin{Proposition} \label{P#1}}
\newcommand{\eP}{\end{Proposition}}

\newcommand{\bLemma}[1]{
\begin{Lemma} \label{L#1} }
\newcommand{\eL}{\end{Lemma}}

\newcommand{\bCorollary}[1]{
\begin{Corollary} \label{C#1} }
\newcommand{\eC}{\end{Corollary}}

\newcommand{\bRemark}[1]{
\begin{Remark} \label{R#1} }
\newcommand{\eR}{\end{Remark}}

\newcommand{\bDefinition}[1]{
\begin{Definition} \label{D#1} }
\newcommand{\eD}{\end{Definition}}

\newcommand{\Del}{\Delta_x}

\newcommand{\tv}{{\tilde{\vv}}}
\newcommand{\tP}{ {\tilde{\mathbb{P}}}}
\newcommand{\tQ}{{\tilde{\mathbb{Q}}}}

\newcommand{\data}{{\rm data}}

\newcommand{\tor}{\mathcal{T}^3}

\newcommand{\bFormula}[1]{
\begin{equation} \label{#1}}
\newcommand{\eF}{\end{equation}}

\newcommand{\Ov}[1]{\overline{#1}}

\newcommand{\DC}{C^\infty_c}

\newcommand{\vc}[1]{{\bf #1}}

\newcommand{\Div}{{\rm div}_x}
\newcommand{\Grad}{\nabla_x}

\newcommand{\dx}{\,{\rm d} {x}}

\newcommand{\dt}{\,{\rm d} t }

\newcommand{\intO}[1]{\int_{\mathcal{T}^3} #1 \ \dx}

\newcommand{\vv}{\vc{v}}

\newcommand{\D}{{\rm d}}

\definecolor{Cgrey}{rgb}{0.85,0.85,0.85}
\definecolor{Cblue}{rgb}{0.50,0.85,0.85}
\definecolor{Cred}{rgb}{1,0,0}
\definecolor{fancy}{rgb}{0.10,0.85,0.10}

\newcommand\Cbox[2]{%
    \newbox\contentbox%
    \newbox\bkgdbox%
    \setbox\contentbox\hbox to \hsize{%
        \vtop{
            \kern\columnsep
            \hbox to \hsize{%
                \kern\columnsep%
                \advance\hsize by -2\columnsep%
                \setlength{\textwidth}{\hsize}%
                \vbox{
                    \parskip=\baselineskip
                    \parindent=0bp
                    #2
                }%
                \kern\columnsep%
            }%
            \kern\columnsep%
        }%
    }%
    \setbox\bkgdbox\vbox{
        \color{#1}
        \hrule width  \wd\contentbox %
               height \ht\contentbox %
               depth  \dp\contentbox
        \color{black}
    }%
    \wd\bkgdbox=0bp%
    \vbox{\hbox to \hsize{\box\bkgdbox\box\contentbox}}%
    \vskip\baselineskip%
}

\date{}

\newcounter{marnote}

\begin{document}


\title{On a hyperbolic system arising in liquid crystals modeling}

\author{Eduard Feireisl\thanks{Institute of Mathematics of the Academy of Sciences of the Czech Republic, \v Zitn\' a 25, 
    CZ-115 67 Praha 1, Czech Republic {\it (feireisl@math.cas.cz)}},\hspace{-1ex}
\quad Elisabetta Rocca\thanks{Universit\`a degli Studi di Pavia, Dipartimento di Matematica, 
   Via Ferrata 5, 27100, Pavia, Italy\hspace{15ex} {\it (elisabetta.rocca@unipv.it)}}\\
 Giulio Schimperna\thanks{Universit\`a degli Studi di Pavia, Dipartimento di Matematica, 
  Via Ferrata 5, 27100, Pavia, Italy {\it (giusch04@unipv.it)}},\hspace{-1ex} \quad 
   Arghir Zarnescu\thanks{ IKERBASQUE, Basque Foundation for Science, Maria Diaz de Haro 3,
   48013, Bilbao, Bizkaia, Spain}\,\,\thanks{BCAM, Basque Center for Applied Mathematics, Mazarredo 14, 
 E48009 Bilbao, Bizkaia, Spain {\it (azarnescu@bcamath.org)}}\,\,
   \thanks{ ``Simion Stoilow" Institute of the Romanian Academy, 21 Calea Grivi\c{t}ei,  010702 Bucharest, Romania}}

\date{\today}

\maketitle

\begin{abstract}

We consider a model of liquid crystals, based on a nonlinear hyperbolic system of differential equations, that 
 represents an inviscid version of the model proposed by Qian and Sheng. A new concept of
\emph{dissipative solution} is proposed, for which a global-in-time existence theorem is shown. The dissipative solutions enjoy the following properties:

(i) they exist globally in time for any finite energy initial data;

(ii) dissipative solutions enjoying certain smoothness are classical solutions;

(iii) a dissipative solution coincides with a strong solution originating from the same initial data as long as the latter exists.

\end{abstract}



\section{Introduction}
\label{i}

In this article we study a system modelling the hydrodynamics of nematic liquid crystals in the Q-tensor framework. The system is an inviscid version 
of the equations proposed by T.~Qian and P.~Sheng in \cite{QiaShe} and studied analytically in \cite{FaZa}. It was proposed by 
F.~Gay-Balmaz and C.~Tronci in \cite{GaTr} as a simplification of the Qian-Sheng model, that captures its 
essential features and exhibits a number of interesting conservation and geometric properties that could be
relevant in particular to describing  defect patterns, thanks to the (presumptively) singular character of the equations.

The most characteristic specific feature of this model is the presence of an inertial term that appears as a second-order material derivative.
This term provides a hyperbolic character to the equations and  is the main source of difficulties in the analysis. It should be noted that this 
second material derivative also appears in the more commonly used Ericken-Leslie model of liquid crystals, but there it generates even worse effects due to the additional presence of the unit-length constraint. 

The system we are studying is:
\begin{align}
\label{i1}
&\Div \vv = 0 \\
\label{i2} 
&\partial_t \vv + \vv \cdot \Grad \vv + \Grad \Pi = - \Div \left(  \Grad \mathbb{Q} \odot \Grad \mathbb{Q} \right)   \\
\label{i3}
&\partial_t \mathbb{Q} + \vv \cdot \Grad \mathbb{Q} = \mathbb{P} \\
\label{i4}
&\partial_t \mathbb{P} + \vv \cdot \Grad \mathbb{P} = - \frac{\partial \mathcal{F} }{\partial \mathbb{Q} } + \Delta \mathbb{Q} - \lambda \mathbb{I}
\end{align}
For the sake of simplicity, we restrict ourselves to the periodic boundary conditions, 
for which the underlying spatial domain may be identified with the flat torus:
\begin{equation} \label{i5}
\tor = \left( [-\pi, \pi]|_{\{ - \pi, \pi \}} \right)^3.
\end{equation}

The system (\ref{i1}), (\ref{i2}) is the standard Euler system for the fluid velocity $\vv = \vv(t,x) \in R^3$, coupled via a nonlinear forcing term with a wave-like
equation (\ref{i3}), (\ref{i4}) governing the time evolution of the $Q-$tensor $\mathbb{Q}  = \mathbb{Q}(t,x) \in R^{3\times 3}_{0, {\rm sym}}$ - a symmetric traceless matrix.
The pressure $\Pi$ and the factor $\lambda \mathbb{I}$ may be seen as Lagrange multipliers compensating the deviation of the motion from the 
divergenceless and zero-trace state, respectively. The problem admits a natural energy functional
\begin{equation} \label{i6}
\mathcal{E}[\vv, \mathbb{P}, \mathbb{Q} ] = \intO{ \left[ \frac{1}{2} |\vv|^2 + \frac{1}{2} |\mathbb{P}|^2 + \frac{1}{2} |\Grad \mathbb{Q}|^2 + \mathcal{F}(\mathbb{Q}) \right] }.
\end{equation}
It is easy to check, by multiplying (\ref{i2}) on $\vv$, (\ref{i4}) on $\mathbb{P}$, and integrating the resulting sum by parts, that
the total energy is conserved
\begin{equation} \label{i7}
\frac{{\rm d}}{{\rm d}t} \mathcal{E}[\vv, \mathbb{P}, \mathbb{Q}] = 0
\end{equation}
for any smooth solution of (\ref{i1}--\ref{i5}).

\bigskip
One of the main interesting qualitative features of the whole system, \eqref{i1}-\eqref{i4} is that it can be regarded as an 
{\it ``extended Euler system"} and this was the main motivation for proposing it in \cite{GaTr}. Indeed, as pointed in \cite{GaTr}
if one defines the {\it ``extended circulation vector"} \footnote{ here and in the following we assume Einstein summation convention, of summation over repeated indices}
$$\mathcal{C}_{QS}:=\vv+\mathbb{P}_{ij}\Grad \mathbb{Q}_{ij}$$ and the {\it ``extended vorticity"}
$$\bar \omega_{QS}:=\Grad\times \mathcal{C}_{QS}$$ then we have the Euler-like equation 
\begin{equation}\label{extEuler}
\partial_t \bar \omega_{QS}+\Grad\times(\vv \times\bar\omega_{QS})=0
\end{equation}
In addition there hold a {\it ``circulation theorem"} 
$$\frac{{\rm d}}{{\rm d}t}\oint_{\Gamma(t)} \vv+\mathbb{P}_{ij}\Grad\mathbb{Q}_{ij} \,d\sigma=0$$
where $\Gamma(t)$ is a closed path moving with velocity $v$ 
and the {\it ``helicity conservation" } 
$$\frac{{\rm d}}{{\rm d}t}\intO{ \mathbb{C}_{QS}\cdot\bar\omega_{QS}}=0$$ 
as a direct consequence of the equation \eqref{extEuler}. 
Moreover the equations allow for the existence of singular vortex structures (see \cite{GaTr} for details).

\bigskip
Our goal in this paper is to focus on quantitative aspects and  study the existence of \emph{global-in-time} solutions to problem (\ref{i1}--\ref{i5}).
This may seem a rather ambitious task as the problem
is highly non-linear involving the incompressible Euler system for which the existence of \emph{physically admissible} solutions is an open problem
even in the class of weak solutions, see however Wiedemann \cite{Wied}, Sz\' ekelyhidi--Wiedemann \cite{SzeWie}. To circumvent these well-known 
difficulties, we introduce a new class of \emph{dissipative solutions} inspired by a similar concept introduced by
P.L.~Lions \cite{LI4a} in the context of Euler flow. The dissipative solutions enjoy the following properties:
\begin{itemize}
\item Any classical solution of problem (\ref{i1}--\ref{i5}) is a dissipative solution.
\item Any (sufficiently) smooth dissipative solution is a classical solution.
\item A dissipative solution coincides with the classical solution emanating from the same initial data as long as the latter exists.
\end{itemize}

Our strategy is based on the following steps. First, we establish existence of smooth solutions of (\ref{i1}--\ref{i5})
defined on a possibly short time interval the length of which depends on the norm of the initial data in certain Sobolev spaces. This will be done in an entirely
standard way by the energy method well developed in the theory of hyperbolic conservation laws. Next, we introduce 
the concept of dissipative solution. Roughly speaking, the dissipative solutions will satisfy the system of equations (\ref{i1}--\ref{i4}) in the sense of distributions, where the right-hand side
of (\ref{i1}), (\ref{i4}) will contain two extra terms playing the role of \emph{defect measures}. In addition, the dissipative solutions will satisfy the
\emph{energy inequality}
\[
\mathcal{E}[\vv, \mathbb{P}, \mathbb{Q}](\tau) + \mathcal{D}(\tau) \leq  \mathcal{E}[\vv, \mathbb{P}, \mathbb{Q}](0)\ \mbox{for a.a.}\ \tau > 0,
\]
with a \emph{dissipation defect} $\mathcal{D}$ dominating, in a certain sense specified below, the defect measure in the equations. Next, we derive a relative
energy inequality playing the role of a ``distance'' between a dissipative solution and any sufficiently smooth process. Similarly to \cite{FeNoSun}, we use
the relative energy to show the weak-strong uniqueness property for the class of dissipative solutions. Finally, we observe that the same procedure
used in the construction of local smooth solutions gives rise to a dissipative solution.

\bigskip
The  main results and the organization of the paper are as follows:
\begin{itemize}
\item In Section \ref{E}, we establish the existence of classical solutions to the initial-value problem
for (\ref{i1}--\ref{i5}) on a possibly short time interval, see Theorem \ref{TE1}.
\item  In Section \ref{LSD}, we introduce the concept of dissipative solution
and show their global-in-time existence in Theorem \ref{TLL1}
\item  Finally, in Section \ref{L}, we show the weak--strong uniqueness principle in the class of dissipative solutions in Theorem \ref{TWS}.
\end{itemize}

\smallskip

{\bf Notations and conventions:} For $A,B$ two $3\times 3$ matrices we denote the inner product on the space of matrices as $A:B=\textrm{tr}(AB)$.  The product $\nabla A\otimes \nabla B$ is a matrix with $ij$ component $\partial_i A: \partial_j B$. 

If $M(x)$ is a $3\times 3$-matrix, then $|M|$ denotes the Frobenius norm of the matrix, i.e. $|M|=\sqrt{M:M^t}$. Furthermore $\nabla \cdot M$ stands for the vector field $(\sum_{j=1}^3 \frac{\partial M _{ij}}{\partial x_j})_{i=1\dots 3}$. 

If $\vc{v}$ is a $3$-dimensional vector and $\mathbb{Q}$ is a $3\times 3$ matrix then $\vc{v}\otimes\mathbb{Q}$ is a third order tensor with components $v_iQ_{kl}$ with $i,k,l\in \{1,2,3\}$. Also if $\mathbb{P}$ is a  $3\times 3$ then $\Grad \mathbb{P}$ is a third-order tensor and we denote by $v\otimes \mathbb{Q}:\Grad\mathbb{P}$ the scalar $\sum_{i,j,k=1}^3 v_k Q_{ij} \frac{\partial P_{ij}}{\partial  x_k}$.

 Occasionally for simplicity we  will write $\partial\mathcal{G}(\mathbb{Q})$ for $\frac{\partial \mathcal{G}}{\partial \mathbb{Q}}(\mathbb{Q})$.


\section{The existence of local-in-time smooth solutions}
\label{E}

We establish the existence of local-in-time solutions in the Sobolev framework $W^{s,2}(\mathcal{T}^N)$ of functions with derivatives
up to order $s$ square integrable in $\mathcal{T}^N$. Let us start by recalling the following nowadays standard results,
see e.g. Majda \cite[Proposition 2.1]{Majd}:
\begin{enumerate}
\item For $u,v \in W^{s,2} \cap L^\infty(\tor)$ and $\alpha$ a multi-index with $|\alpha| \leq s$
\begin{equation}\label{E5}
\left\| \partial^\alpha_x (u v) \right\|_{L^2(\tor)} \leq c_s \left( \| u \|_{L^\infty(\tor)} \| \nabla_x^s v \|_{L^{2}(\tor)} +
\| v \|_{L^\infty(\tor)} \|  \nabla_x^s u \|_{L^{2}(\tor)} \right).
\end{equation}
Here and below $ \nabla_x^s v$ denotes the tensor of the partial derivatives of $v$ of order equal to $s$.

\item For $u \in W^{s,2}(\tor)$, $\nabla_x u \in L^\infty(\tor)$, $v \in W^{s-1,2} \cap L^\infty (\tor)$ and $|\alpha| \leq s$
\begin{equation}\label{E6}
\left\| \partial^\alpha_x (uv) - u \partial^\alpha_x v \right\|_{L^2(\tor)} \leq c_s
\left( \| \nabla_x u \|_{L^\infty(\tor)} \| \nabla^{s-1}_x v \|_{L^2(\tor)} +
 \| \nabla^s_x u \|_{L^2(\tor)} \| v \|_{L^\infty(\tor)} \right).
\end{equation}

\item For $u \in W^{s,2} \cap C(\tor)$, and $F$ $s$-times continuously differentiable function on an open neighborhood of the compact set $G =
{\rm range}[u]$, $1\le |\alpha| \leq s$,
\begin{equation} \label{E7}
\left\| \partial^\alpha_x F(u) \right\|_{L^2(\tor)} \leq c_s \| \partial_u F \|_{C^{s-1}(G)} \| u \|^{|\alpha| - 1}_{L^\infty(\tor)} \|
\partial^\alpha_x u \|_{L^2(\tor)}.
\end{equation}

\end{enumerate}

We will need  to make a number of assumptions about the potential $\mathcal{F}$, namely:
\begin{itemize}
\item[{\bf (A1)}] 
The function  $\mathcal{F}:R^{3\times 3}_{0, {\rm sym}}\to R$ is {\it isotropic}, i.e.
\[
 \mathcal{F}(\mathbb{Q})=\mathcal{F}(\mathbb{R}\mathbb{Q}\mathbb{R}^t),~\forall \mathbb{R}\in O(3).
\]
\item[{\bf (A2)}] 
There exists $\Lambda\ge 0$ such that $\mathcal{F}+\Lambda|\mathbb{Q}|^2$ is a strictly convex non-negative function.
Hence we will set
\[
  \mathcal{G}(\mathbb{Q}) := \mathcal{F}(\mathbb{Q})+\Lambda|\mathbb{Q}|^2.
\]
\item[{\bf (A3)}] 
There exists a constant $\bar C>0$ such that 
\begin{equation}\label{assumptionA3}
 |\partial \mathcal{F}(\mathbb{Q})|\le \bar C(1+|\mathbb{Q}|^q), 
\end{equation}
for some exponent $q<5$.
\end{itemize}
\begin{Remark}\label{growth}
  The first assumption (A1) is motivated by physical reasons. Assumption (A2),
  usually termed as ``$\Lambda$-convexity'' is also natural and allows the occurrence
  of potentials with many wells. Regarding (A3), this hypotheses is taken just for
  convenience in order to simplify the subsequent computations. Indeed, refining a bit 
  the proof one could treat also potentials with faster (polynomial) growth
  at infinity, the only complication being the possible appearence of a 
  further defect term in the energy balance (cf.~\eqref{WSE} below). Indeed,
  the function $G(\mathbb{Q})$ can be controlled only in $L^1$ in that case.
  On the other hand, managing potentials with exponential growth at infinity 
  or ``singular'' potentials (i.e., functions $\mathcal{F}$ being identically
  infinity outside a bounded set, like for instance the
  Ball-Majumdar potential considered in \cite{BM,FRSZ}) may be more delicate
  because in that situation one would also face the occurrence of a 
  further measure-valued term in \eqref{i4} resulting as one takes the
  limit of $\partial \mathcal{F}(\mathbb{Q})$.
  We finally observe that an example of function that satisfies (A1)-(A3)
  and is often used is given by
  $$
    \mathcal{F}(\mathbb{Q})=\frac{a}{2}|\mathbb{Q}|^2+\frac{b}{3}\textrm{tr}(\mathbb{Q}^3)+\frac{c}{4}|\mathbb{Q}|^4
  $$
  where $a,b\in R$ and $c>0$.
\end{Remark}
Our goal is to show the following result:

\Cbox{Cgrey}{
\begin{Theorem} \label{TE1}

Let $s\ge 3$ and $\mathcal{F} \in C^{s + 1}(R^{3 \times 3}_{0,{\rm sym}}; R)$ be given and satisfying the assumptions $(A1)$ through $(A3)$ before.
Consider the initial data
\begin{equation} \label{i5a}
\vc{v}(0, \cdot) = \vc{v}_0 \in W^{s,2}(\tor; R^3), \ \mathbb{P}(0, \cdot) = \mathbb{P}_0 \in
W^{s,2}(\tor; R^{3 \times 3}_{0,{\rm sym}}),\ \mathbb{Q}(0, \cdot) = \mathbb{Q}_0 \in
W^{s + 1,2}(\tor; R^{3 \times 3}_{0,{\rm sym}})
\end{equation}
such that
\[
\Div \vc{v}_0  = 0.
\]
Then there exists $T_0 > 0$ depending solely on the norm of the initial data in the aforementioned spaces such that problem (\ref{i1}--\ref{i4}) admits
a strong solution in $[0, T_0] \times \tor$, unique in the class
\[
\begin{split}
\vc{v} &\in C([0,T_0]; W^{s,2}(\tor; R^3), \\
\mathbb{P} &\in C([0,T_0]; W^{s,2}(\tor; R^{3 \times 3}_{0,{\rm sym}}) ),\ \mathbb{Q} \in C([0,T_0]; W^{s + 1,2}(\tor; R^{3 \times 3}_{0,{\rm sym}}) ).
\end{split}
\]
\end{Theorem}
}
The rest of this section is devoted to the proof of Theorem \ref{TE1}.


\subsection{A priori bounds}

We start by establishing {\it a priori} estimates for smooth solutions of problem (\ref{i1}--\ref{i5}), (\ref{i5a}).


\subsubsection{Energy estimates}


First of all, take the scalar product of (\ref{i4}) with $\mathbb{P}$ and use (\ref{i1}), (\ref{i3}) obtaining
\begin{equation} \label{E7b}
\frac{\D }{\dt} \intO{ \left[ \frac{1}{2} |\mathbb{P}|^2 + \frac{1}{2} |\Grad \mathbb{Q}|^2 + \mathcal{G}(\mathbb{Q}) \right] } =
\intO{ (\vc{v} \cdot \Grad \mathbb{Q}) : \Delta \mathbb{Q} } + 2 \Lambda \intO{ \mathbb{Q} : \mathbb{P} }.
\end{equation}

Next, take the scalar product of (\ref{i2}) with $\vc{v}$, use \eqref{i1} and add the resulting expression to (\ref{E7b}) to get
\begin{equation} \label{E7c}
\frac{\D }{\dt} E[\vv, \mathbb{P}, \mathbb{Q}] =
 2 \Lambda \intO{ \mathbb{Q} : \mathbb{P} },
\end{equation}
where we have set
\begin{equation} \label{defi:E}
 E[\vv, \mathbb{P}, \mathbb{Q}] := \intO{ \left[ \frac{1}{2} |\vc{v}|^2 + \frac{1}{2} |\mathbb{P}|^2 + \frac{1}{2} |\Grad \mathbb{Q}|^2 + \mathcal{G}(\mathbb{Q}) \right] } .
\end{equation}

Consequently, applying Gronwall's lemma and recalling that $\mathcal{G}$ was assumed to be strictly
convex (cf.~(A2)), we deduce the energy bounds
\begin{equation} \label{E7d}
\begin{split}
\sup_{t \in [0,T]} \| \vc{v}(t, \cdot) \|_{L^2(\tor; R^3)} &\leq c(T, \data), \\
\sup_{t \in [0,T]} \| \mathbb{P}(t, \cdot) \|_{L^2(\tor; R^{3 \times 3})} &\leq c(T, \data), \\
\sup_{t \in [0,T]} \| \mathbb{Q}(t, \cdot) \|_{W^{1,2}(\tor; R^{3 \times 3})} &\leq c(T, \data).
\end{split}
\end{equation}
\begin{Remark} \label{Rdd1}
The energy estimates (\ref{E7d}) are uniform on any bounded time interval $(0,T)$.
\end{Remark}


\subsubsection{Higher order estimates}
\label{Eh}

Rewrite (\ref{i2}) as
\begin{equation} \label{A1}
\partial_t \vc{v} + \vc{v} \cdot \Grad \vc{v} + \Grad \Pi = - \Del \mathbb{Q} : \Grad \mathbb{Q} - \frac12 \Grad |\Grad \mathbb{Q}|^2
\end{equation}
and apply $\partial^\alpha_x$, $1\le |\alpha| \leq s$, to both sides of (\ref{A1}) to obtain
\[
\begin{split}
\partial_t \left( \partial^\alpha_x \vc{v} \right) &+ \vc{v} \cdot \Grad \left( \partial^\alpha_x \vc{v} \right)  + \Grad (\partial^\alpha_x \Pi) =
- \left( \partial^\alpha_x \Del \mathbb{Q} \right): \Grad \mathbb{Q} - - \frac12 \Grad \partial^\alpha_x |\Grad \mathbb{Q}|^2\\
&+ \vc{v}  \cdot \left( \Grad\partial^\alpha_x \vc{v} \right) - \partial^\alpha_x \left( \vc{v} \cdot \Grad \vc{v} \right)
 - \partial^\alpha_x \left( \Del \mathbb{Q} : \Grad \mathbb{Q} \right) + \left( \partial^\alpha_x \Del \mathbb{Q} \right): \Grad \mathbb{Q}.
\end{split}
\]
Next,
take the scalar product of this equation with $\partial^\alpha_x \vc{v}$ and integrate over $\tor$ using (\ref{i1}) to obtain
\begin{equation} \label{A2}
\begin{split}
&\frac{1}{2} \frac{{\rm d}}{{\rm d}t} \intO{ | \partial^\alpha_x \vc{v} |^2 } =  - \intO{ \left( \partial^\alpha_x \Del \mathbb{Q} \right): \Grad \mathbb{Q}
\cdot \left( \partial^\alpha_x \vc{v} \right) } \\
& - \intO{ \left[ \partial^\alpha_x \left( \vc{v} \cdot \Grad \vc{v} \right) - \vc{v}  \cdot \left( \partial^\alpha_x \Grad \vc{v} \right) \right]
\cdot \partial^\alpha_x \vc{v} } \\
- &\intO{ \left[ \partial^\alpha_x \left( \Del \mathbb{Q} : \Grad \mathbb{Q} \right) -  \left( \partial^\alpha_x \Del \mathbb{Q} \right): \Grad \mathbb{Q}
\right] \cdot \partial^\alpha_x \vc{v} },
\end{split}
\end{equation}
where, in accordance with (\ref{E6}),
\begin{equation} \label{A3}
\begin{split}
\left| \intO{ \left[ \partial^\alpha_x \left( \vc{v} \cdot \Grad \vc{v} \right) - \vc{v}  \cdot \left( \partial^\alpha_x \Grad \vc{v} \right) \right]
\cdot \partial^\alpha_x \vc{v} } \right|
\leq c  \left\| \Grad \vc{v} \right\|_{L^\infty(\tor)} \left\| \Grad^s \vc{v} \right\|_{L^2(\tor)}
\left\| \partial^\alpha_x \vc{v} \right\|_{L^2(\tor)},
\end{split}
\end{equation}
and, by the same token
\begin{equation} \label{A4}
\begin{split}
&\left| \intO{ \left[ \partial^\alpha_x \left( \Del \mathbb{Q} : \Grad \mathbb{Q} \right) -  \left( \partial^\alpha_x \Del \mathbb{Q} \right): \Grad \mathbb{Q}
\right] \cdot \partial^\alpha_x \vc{v} } \right| \leq \\ &\leq \left[ \left\| \Grad^2 \mathbb{Q} \right\|_{L^\infty(\tor)} \left\|
\Grad^{s+1} \mathbb{Q} \right\|_{L^2(\tor)} + \left\| \Del \mathbb{Q} \right\|_{L^\infty(\tor)} \left\|
\Grad^{s+1} \mathbb{Q} \right\|_{L^2(\tor)}\right] \left\| \partial^\alpha_x \vc{v} \right\|_{L^2(\tor)}.
\end{split}
\end{equation}
The next step is to apply $\partial^\alpha_x$ to equation (\ref{i4}):
\[
\partial_t \left( \partial^\alpha_x \mathbb{P} \right) + \vc{v} \cdot \Grad \left( \partial^\alpha_x \mathbb{P} \right) = \partial^\alpha _x ( \Del  \mathbb{Q} ) -\partial^\alpha_x \left( \partial \mathcal{F} (\mathbb{Q}) \right)  +
\vc{v} \cdot \Grad \left( \partial^\alpha_x \mathbb{P} \right) - \partial^\alpha_x \left( \vc{v} \cdot \Grad   \mathbb{P} \right).
\]
Now, similarly to the above, take the scalar product with $\partial^\alpha_x \mathbb{P}$ and integrate by parts to obtain
\begin{equation} \label{A5}
\begin{split}
\frac{1}{2} \frac{{\rm d}}{{\rm d}t} \intO{ | \partial^\alpha_x \mathbb{P} |^2 } & =
\intO{ \partial^\alpha _x ( \Del  \mathbb{Q} ): \partial^\alpha_x \mathbb{P}}  - \intO{ \partial^\alpha_x \left( \partial \mathcal{F} (\mathbb{Q}) \right) : \partial^\alpha_x \mathbb{P} } \\ + & \intO{\left( \vc{v} \cdot \Grad \left( \partial^\alpha_x \mathbb{P} \right) - \partial^\alpha_x \left( \vc{v} \cdot \Grad   \mathbb{P} \right)\right): \partial^\alpha_x \mathbb{P} },
\end{split}
\end{equation}
where, by virtue of (\ref{E6}),
\[
\begin{split}
&\left| \intO{ \vc{v} \cdot \Grad \left( \partial^\alpha_x \mathbb{P} \right) - \partial^\alpha_x \left( \vc{v} \cdot \Grad   \mathbb{P} \right): \partial^\alpha_x \mathbb{P} } \right| \\
&\leq c \left[ \left\| \Grad \vc{v} \right\|_{L^\infty(\tor)} \left\| \Grad^s \mathbb{P} \right\|_{L^2(\tor)} + \left\| \Grad \mathbb{P} \right\|_{L^\infty(\tor)} \left\| \Grad^s \vc{v} \right\|_{L^2(\tor)}
 \right] \left\| \partial^\alpha_x \mathbb{P} \right\|_{L^2(\tor)}.
\end{split}
\]
Next, expressing $\mathbb{P}$ by means of equation (\ref{i3}), we get
\begin{equation} \label{A6}
\begin{split}
&\intO{ \partial^\alpha _x ( \Del  \mathbb{Q} ): \partial^\alpha_x \mathbb{P}} = \intO{ \partial^\alpha _x ( \Del  \mathbb{Q} ): \partial^\alpha_x
\partial_t \mathbb{Q} } + \intO{ \partial^\alpha _x ( \Del  \mathbb{Q} ): \partial^\alpha_x
\left( \vc{v} \cdot \Grad \mathbb{Q} \right) }\\
&= - \frac{1}{2}\frac{{\rm d}}{{\rm d}t} \intO{ | \Grad \partial^\alpha_x \mathbb{Q} |^2 }
+ \intO{ \partial^\alpha _x ( \Del  \mathbb{Q} ): \Grad \mathbb{Q} \cdot \partial^\alpha_x \vc{v} } +
\intO{ \Del \left( \partial^\alpha _x  \mathbb{Q} \right): \Grad \left( \partial^\alpha_x \mathbb{Q} \right) \cdot \vc{v} }\\
&+ \intO{ \partial^\alpha _x ( \Del  \mathbb{Q} ): \left[ \partial^\alpha_x
\left( \vc{v} \cdot \Grad \mathbb{Q} \right) - \Grad \mathbb{Q} \cdot \partial^\alpha_x \vc{v} -  \partial^\alpha_x \left( \Grad \mathbb{Q} \right) \cdot \vc{v}
 \right] }.
\end{split}
\end{equation}
Note carefully that
\begin{itemize}
\item the integral
\[
\intO{ \partial^\alpha _x ( \Del  \mathbb{Q} ): \Grad \mathbb{Q} \cdot \partial^\alpha_x \vc{v} }
\]
cancels out with its counterpart in (\ref{A2});
\item $\vc{v}$ is solenoidal, whence
\[
\begin{split}
\intO{ \Del \left( \partial^\alpha _x  \mathbb{Q} \right): \Grad \left( \partial^\alpha_x \mathbb{Q} \right) \cdot \vc{v} } &=
\intO{ \Div \left[ \Grad \left( \partial^\alpha_x \mathbb{Q} \right) \odot \Grad \left( \partial^\alpha_x \mathbb{Q} \right) \right]
\cdot \vc{v} }\\
&= - \intO{ \left[ \Grad \left( \partial^\alpha_x \mathbb{Q} \right) \odot \Grad \left( \partial^\alpha_x \mathbb{Q} \right) \right]
: \Grad \vc{v} }\\
&  \le \|\Grad \mathbb{Q} \|_{W^{s,2}(\tor)}^2\|\Grad v\|_{L^\infty(\tor)}\le   \|\Grad \mathbb{Q}\|_{W^{s,2}(\tor)}^2\|v\|_{W^{s,2}(\tor)}
\end{split}
\] for $s\ge 3$.
\end{itemize}
Thus it remains to handle
\[
\begin{split}
&\intO{ \partial^\alpha _x ( \Del  \mathbb{Q} ): \left[ \partial^\alpha_x
\left( \vc{v} \cdot \Grad \mathbb{Q} \right) - \Grad \mathbb{Q} \cdot \partial^\alpha_x \vc{v} -  \partial^\alpha_x \left( \Grad \mathbb{Q} \right) \cdot \vc{v}
 \right] }\\
& =- \intO{ \partial^\alpha _x ( \Grad  \mathbb{Q} ): \Grad \left[ \partial^\alpha_x
\left( \vc{v} \cdot \Grad \mathbb{Q} \right) - \Grad \mathbb{Q} \cdot \partial^\alpha_x \vc{v} -  \partial^\alpha_x \left( \Grad \mathbb{Q} \right) \cdot \vc{v}
 \right] }.
\end{split}
\]
In view of (\ref{A6}), it is enough to control the norm
\[
\left\| \Grad \left[ \partial^\alpha_x
\left( \vc{v} \cdot \Grad \mathbb{Q} \right) - \Grad \mathbb{Q} \cdot \partial^\alpha_x \vc{v} -  \partial^\alpha_x \left( \Grad \mathbb{Q} \right) \cdot \vc{v}
 \right]
 \right\|_{L^2(\tor)}.
\]
We have
\[
\partial^\alpha_x
\left( \vc{v} \cdot \Grad \mathbb{Q} \right) - \Grad \mathbb{Q} \cdot \partial^\alpha_x \vc{v} -  \partial^\alpha_x \left( \Grad \mathbb{Q} \right) \cdot \vc{v}
= \sum \partial^{\beta_1}_x \vc{v} \cdot \partial^{\beta_2}_x ( \Grad \mathbb{Q} ),
\]
where the sum is extended to all the possible couples $\beta_1$, $\beta_2$ satisfying
\begin{equation} \label{A7}
|\beta_1| \geq 1,\ |\beta_2 | \geq 1, \ \beta_1 + \beta_2 = \alpha.
\end{equation}
Hence,
\begin{equation} \label{A8}
\Grad \left[ \partial^\alpha_x
\left( \vc{v} \cdot \Grad \mathbb{Q} \right) - \Grad \mathbb{Q} \cdot \partial^\alpha_x \vc{v} -  \partial^\alpha_x \left( \Grad \mathbb{Q} \right) \cdot \vc{v}
\right]
= \sum \partial^{\beta_1 + 1}_x \vc{v} \partial^{\beta_2}_x ( \Grad \mathbb{Q} ) +
\sum \partial^{\beta_1}_x \vc{v} \partial^{\beta_2+1}_x ( \Grad \mathbb{Q} ),
\end{equation}
where with some abuse of notation we have noted, for instance,
as $\partial^{\beta_1 + 1}_x \vc{v}$ the gradient of $\partial^{\beta_1}_x \vc{v}$. 
In particular, it turns out that $|\beta_i|+1\le s$.

Let us just focus on the first type of term, $\partial^{\beta_1 + 1}_x \vc{v} \partial^{\beta_2}_x ( \Grad \mathbb{Q} ) $, and consider two cases:
\begin{itemize}
\item If $|\beta_1|+1=s$ and $|\beta_2|=1$ then: 
\[
\|\partial_x^{\beta_1+1} \vc{v} \partial_x^{\beta_2}\Grad \mathbb{Q}\|_{L^2(\tor)}\le \|\vc{v}\|_{W^{s,2}(\tor)}
   \|\Grad\Grad\mathbb{Q}\|_{L^\infty(\tor)}\le  \|\vc{v}\|_{W^{s,2}(\tor)}\|\Grad\mathbb{Q}\|_{W^{s,2}(\tor)},
\] 
provided that $s\ge 3$.
\item If $2\le |\beta_1|+1\le s-1$, $2\le |\beta_2|\le s-1$ then:
\[
  \|\partial_x^{\beta_1+x} \vc{v}\partial_x^{\beta_2}\Grad\mathbb{Q}\|_{L^2(\tor)}
    \le \|\partial_x^{\beta_1+x} \vc{v} \|_{L^4(\tor)}\|\partial_x^{\beta_2}\Grad\mathbb{Q}\|_{L^4(\tor)}.
\]
On the other hand we have out of the Gagliardo-Nirenberg inequalities:
\[
  \|\partial^\gamma_x f\|_{L^4(\tor)}\le \|f\|_{L^2}^{1-a}\|\partial^s_x f\|_{L^2}
   + \|f\|_{L^2}
\] 
provided that $\frac{1}{4}=\frac{|\gamma|}{3}+a(\frac 12 -\frac s3)+\frac{(1-a)}{2}$ and $\frac{|\gamma|}s\le a\le 1$
which holds in our case for $2\le |\gamma|\le s-1$.
\end{itemize}
The above considerations suffice also for bounding the second type of term in \eqref{A8},
namely
$\partial^{\beta_1}_x \vc{v} \cdot \partial^{\beta_2+1}_x ( \Grad \mathbb{Q} )$
just by interswitching $v$ with $\Grad\mathbb{Q}$ and $\beta_1$ with $\beta_2$.

Under these circumstances, we may sum up the relations (\ref{A2}), (\ref{A5}), and (\ref{A6}) to conclude that
\begin{equation} \label{A10}
\begin{split}
&\frac{{\rm d}}{{\rm d}t}  \left[ \left\| \vc{v} \right\|^2_{W^{s,2}(\tor)} +
\left\| \mathbb{P} \right\|_{W^{s,2}(\tor)}^2 + \left\| \mathbb{Q} \right\|_{W^{s+1,2}(\tor)}^2 \right]  \\ \leq
&c \left( \sum_{|\alpha| \leq s } \left| \intO{ \partial^\alpha_x \left( \partial \mathcal{F} (\mathbb{Q}) \right) : \partial^\alpha_x \mathbb{P} } \right|
 + \left[1+ \left\| \vc{v} \right\|^2_{W^{s,2}(\tor)} 
 + \left\| \mathbb{P} \right\|_{W^{s,2}(\tor)}^2 + \left\| \mathbb{Q} \right\|_{W^{s+1,2}(\tor)}^2 \right]^2  \right).
\end{split}
\end{equation}


\subsection{Proof of Theorem \ref{TE1}}

With the bounds established in (\ref{A10}), the proof of Theorem \ref{TE1} can be carried over by means of the Galerkin approximation procedure based on the
trigonometric polynomials similarly to the theory of symmetric hyperbolic systems, cf.~Majda \cite{Majd}. Note that 
all operations performed in Section \ref{Eh} are compatible with such an approximation as the associated
finite-dimensional spaces are closed with respect to all differential operators acting in the
space variable. Accordingly, we may construct a sequence of approximate solutions $[\vv_n, \mathbb{P}_n, \mathbb{Q}_n]$ and perform the limit
for $n \to \infty$ using the bounds (\ref{A10}) on a suitably short time interval $[0, T_0]$.


\section{Dissipative solutions}
\label{LSD}

Following the proof of Theorem \ref{TE1}, we may consider the limit of the approximate solutions $[\vv_n, \mathbb{P}_n, \mathbb{Q}_n]$ 
on an arbitrary time interval $(0,T)$. Indeed, using the a-priori estimates (\ref{E7d}) and comparing terms in \eqref{i2},
\eqref{i3}, \eqref{i4}, we may see that the time derivatives $\partial_t \vc{v}$,
$\partial_t \mathbb{P}$, $\partial_t \mathbb{Q}$ are uniformly bounded in spaces of the 
form $L^p(0,T;X)$ where $p>1$ and $X$ is a suitable Sobolev space of negative order
on the unit torus. Hence, recalling the Assumptions (A2) and (A3) on $\mathcal{F}$, 
from the energy bounds (\ref{E7d}) there follows
\begin{equation} \label{LL1}
\begin{split}
\vv_n &\to \vv \ \mbox{in} \ C_{\rm weak}([0,T]; L^2(\tor; R^3) ),\\
\mathbb{P}_n &\to \mathbb{P} \ \mbox{in}\ C_{\rm weak}([0,T]; L^2(\tor; R^{3 \times 3}) ), \\
\mathbb{Q}_n &\to \mathbb{Q} \ \mbox{in}\ C_{\rm weak}([0,T]; W^{1,2}(\tor; R^{3 \times 3}) )
\ \mbox{and in}\ C([0,T]; L^2(\tor; R^{3 \times 3}) ),
\end{split}
\end{equation}
passing to suitable subsequences as the case may be. To deduce the last property we also used 
the Aubin-Lions lemma. 

Consequently, by virtue of the compact embedding $W^{1,2}(\tor) \hookrightarrow\hookrightarrow L^r$, $1 \leq r < 6$, we also
obtain that there exists $p > 1$ such that for any $\tau \in [0,T]$ there hold:
\begin{equation} \label{LL2}
\begin{split}
\mathcal{G}(\mathbb{Q}_n (\tau, \cdot)) &\to \mathcal{G} (\mathbb{Q} (\tau, \cdot) ) \ \mbox{in}\ L^p(\tor)\\
\partial \mathcal{F}(\mathbb{Q}_n (\tau, \cdot)) &\to \partial \mathcal{F} (\mathbb{Q} (\tau, \cdot) ) \ \mbox{in}\ L^p(\tor).
\end{split}
\end{equation}
Notice that here the growth condition in Assumption~(A3) has also been used 
(see, however, Remark~\ref{growth}).

Next, we observe that
\begin{equation} \label{LL3}
\begin{split}
| \vv_n |^2 &\to \Ov{ |\vv |^2 } \ \mbox{weakly-(*) in}\ L^\infty(0,T; \mathcal{M}^+ (\tor)),\\
| \mathbb{P}_n |^2 &\to \Ov{ |\mathbb{P} |^2 } \ \mbox{weakly-(*) in}\ L^\infty(0,T; \mathcal{M}^+ (\tor)), \\
|\Grad \mathbb{Q}_n |^2 &\to \Ov{ |\Grad \mathbb{Q} |^2 } \ \mbox{weakly-(*) in}\ L^\infty(0,T; \mathcal{M}^+ (\tor)),
\end{split}
\end{equation}
and, accordingly,
\begin{equation} \label{LL4}
\begin{split}
| \vv_n - \vc{v} |^2 &\to \Ov{ |\vv |^2 } - |\vc{v}|^2 \ \mbox{weakly-(*) in}\ L^\infty(0,T; \mathcal{M}^+ (\tor)),\\
| \mathbb{P}_n - \mathbb{P}|^2 &\to \Ov{ |\mathbb{P} |^2 } - |\mathbb{P} |^2 \ \mbox{weakly-(*) in}\ L^\infty(0,T; \mathcal{M}^+ (\tor)), \\
|\Grad \mathbb{Q}_n - \Grad \mathbb{Q} |^2 &\to \Ov{ |\Grad \mathbb{Q} |^2 } - |\Grad \mathbb{Q}|^2 \ \mbox{weakly-(*) in}\ L^\infty(0,T; \mathcal{M}^+ (\tor)).
\end{split}
\end{equation}
Moreover, examining the difference
\begin{equation} \label{LL4b}
v^i_n v^j_n - v^i v^j = (v^i_n - v^i) (v^j_n - v^j ) -  v^i (v^j - v^j_n) - v^j (v^i - v^i_n),
\end{equation}
we find out that
\[
  v^i_n v^j_n - v^i v^j \to \mathcal{R}_{1,1}^{i,j} \ \mbox{weakly-(*) in} \ L^\infty(0,T; \mathcal{M}(\tor) ).
\]
Let now $\zeta\in C(\tor)$ with $\| \zeta \|_\infty \le 1$. Then, testing \eqref{LL4b} by $\zeta$,
we obtain
\begin{equation} \label{LL4c}
  \intO{(v^i_n v^j_n - v^i v^j)\zeta} 
   \le \frac12 \intO { (v^i_n - v^i)^2 |\zeta| } 
   + \frac12 \intO { (v^j_n - v^j)^2 |\zeta| } 
\end{equation}
\[
   - \intO{  \big( v^i (v^j - v^j_n) + v^j (v^i - v^i_n) \big)\zeta}
\]
Hence, letting $n\nearrow\infty$, we obtain
\[
 \intO{\mathcal{R}_{1,1}^{i,j}\zeta} 
  \le \frac12 \lim_{n\nearrow\infty}\intO { (v^i_n - v^i)^2 |\zeta| } 
   + \frac12 \lim_{n\nearrow\infty}\intO { (v^j_n - v^j)^2 |\zeta| },
\] 
where the first integral has in fact to be intended as the integral of the function $\zeta$ 
with respect to the measure $\mathcal{R}_{1,1}^{i,j}$. This convention will be extensively
used also in the sequel.

Hence, passing to the supremum with respect to $\zeta$, recalling \eqref{LL4} and summing
over $i,j$, we arrive at
%
%
%
%
\begin{equation} \label{LL5}
\sum_{i,j} \| \mathcal{R}_{1,1}^{i,j} \|_{\mathcal{M} (\tor) }  \leq 3 \intO{ \left( \Ov{ |\vv|^2  } - |\vv|^2 \right) }.
\end{equation}
Proceeding in a similar way, we can prove that
\[
 \sum_{\alpha,\beta} \left( \partial_i Q^{\alpha,\beta}_n \partial_j Q^{\alpha,\beta}_n 
  - \partial_i Q^{\alpha,\beta} \partial_j Q^{\alpha,\beta} \right) \to \mathcal{R}^{i,j}_{1,2} \ \mbox{weakly-(*) in} \ L^\infty(0,T; \mathcal{M}(\tor) ),
\]
where
\begin{equation} \label{LL6}
  \sum_{i,j}  \| \mathcal{R}_{1,2}^{i,j} \|_{\mathcal{M} (\tor) } \leq c \intO{ \left( \Ov{ |\Grad \mathbb{Q} |^2  } - |\Grad \mathbb{Q}|^2 \right) },
\end{equation}
for some $c>0$, and
\[
v^i_n P^{k,j}_n - v^i P^{k,j} \to \mathcal{R}^{i,j,k}_2 \ \mbox{weakly-(*) in} \ L^\infty(0,T; \mathcal{M}(\tor) ),
\]
where
\begin{equation} \label{LL7}
  \sum_{i,j} \int_0^\tau \| \mathcal{R}_2^{i,j,k} \|_{\mathcal{M} (\tor) } \dt 
    \leq  c \int_0^\tau \intO{ \left( \Ov{ |\vc{v} |^2  } - |\vc{v} |^2 \right) } \dt 
     + c \int_0^\tau \intO{ \left( \Ov{ |\mathbb{P} |^2  } - |\mathbb{P} |^2 \right) } \dt.
\end{equation}
Noting now as $\mathbb{R}^1$ the tensor of measures whose $(i,j)$-entry is 
$\mathcal{R}_{1,1}^{i,j}+\mathcal{R}_{1,2}^{i,j}$
and, respectively, as $\mathbb{R}^2$ the tensor of measures whose $(i,j,k)$-entry is 
$\mathcal{R}_2^{i,j,k}$, we obtain that the 
limit functions $[\vc{v}, \mathbb{P}, \mathbb{Q}]$ satisfy
\begin{equation} \label{WS}
\begin{split}
\intO{ \vv(\tau, \cdot) \cdot \Grad \varphi }  &= 0 \ \mbox{for any}\ \varphi \in \DC(\tor), \\
\left[ \intO{ \vc{v} \cdot \varphi } \right]_{t = 0}^{t = \tau} &= \int_0^\tau \intO{
\vv \cdot \partial_t \varphi + (\vv \otimes \vv) : \Grad \varphi + (\Grad \mathbb{Q} \odot \Grad \mathbb{Q}) : \Grad \varphi } \dt  \\
&+ \int_0^\tau \left< \mathbb{R}^1 ; \Grad \varphi \right> \dt   \\
\mbox{for any}\ \varphi &\in \DC([0,T] \times \tor; R^3), \ \Div \varphi  = 0;\\
\left[ \intO{ \mathbb{Q} : \varphi } \right]_{t = 0}^{t = \tau} &= \int_0^\tau \intO{ \left[
\mathbb{Q} : \partial_t \varphi  + (\vv \otimes \mathbb{Q}) : \Grad \varphi + \mathbb{P}: \varphi \right]} \dt \\
\mbox{for any}\ \varphi &\in \DC([0,T] \times \tor; R^{3 \times 3} ); \\
\left[ \intO{ \mathbb{P} : \varphi} \right]_{t = 0}^{t = \tau} &= \int_0^\tau \intO{ \left[
\mathbb{P}: \partial_t \varphi  + (\vv \otimes \mathbb{P} ): \Grad \varphi  -  \frac{\partial \mathcal{F}(\mathbb{Q}) }{\partial \mathbb{Q} }
: \varphi - \Grad \mathbb{Q} :\Grad \varphi \right] } \dt \\
&+ \int_0^\tau \left< \mathbb{R}^2; \Grad \varphi \right> \dt \\
\mbox{for any} \ \varphi &\in \DC([0,T] \times \tor; R^{3 \times 3}_{0, {\rm sym}}).
\end{split}
\end{equation}
Writing now \eqref{E7c} for the index $n$ and integrating in time over $(0,\tau)$, we infer
\begin{equation} \label{E7c2}
\begin{split}
  E[\vv(\tau), \mathbb{P}(\tau), \mathbb{Q}(\tau)]
  + \Big[ E[\vv_n(\tau), \mathbb{P}_n(\tau), \mathbb{Q}_n(\tau)] -  E[\vv(\tau), \mathbb{P}(\tau), \mathbb{Q}(\tau)] \Big]\\
   = E[\vv_n(0), \mathbb{P}_n(0), \mathbb{Q}_n(0)] + 2 \Lambda \int_0^t \intO{ \mathbb{Q}_n : \mathbb{P}_n } \dt.
\end{split}
\end{equation}
Hence, letting $n\nearrow \infty$, we obtain the 
energy balance
\begin{equation} \label{WSE}
\left[ \intO{ \left[ \frac{1}{2} |\vc{v}|^2 + \frac{1}{2} |\mathbb{P}|^2 + \frac{1}{2} |\Grad \mathbb{Q}|^2 + \mathcal{G}(\mathbb{Q}) \right] }
\right]_{t = 0}^{t = \tau} + \mathcal{D} (\tau) =
2 \Lambda \int_0^\tau \intO{ \mathbb{Q} : \mathbb{P} }\dt
\end{equation}
for any $\tau \in [0,T]$. The function $\mathcal{D} \in L^\infty(0,T)$, $\mathcal{D} \geq 0$ is obtained as the limit of the difference
in square brackets in \eqref{E7c2} and represents a \emph{dissipation defect}. Notice in particular that, owing to \eqref{LL2}, the 
part depending on $\mathcal{G}$ disappears (but cf.~Remark~\ref{growth}). Using then \eqref{LL4} and 
(\ref{LL5}--\ref{LL7}), we have
\begin{equation} \label{DDI}
\int_0^\tau \left[ \| \mathbb{R}^1 (t, \cdot) \|_{\mathcal{M}(\tor)} + \| \mathbb{R}^2 (t, \cdot) \|_{\mathcal{M}(\tor)} \right] \dt
\leq c \int_0^\tau \mathcal{D}(t) \dt, \ \tau \in [0, T].
\end{equation}

The trio of functions $[\vc{v}, \mathbb{P}, \mathbb{Q}]$ belonging to the regularity class (\ref{LL1}) and satisfying
(\ref{WS}--\ref{DDI}) for certain $\mathbb{R}^1$, $\mathbb{R}^2$, $\mathcal{D}$ will be termed \emph{dissipative solution} to problem
(\ref{i1}--\ref{i4}). We have just shown the following result:

\Cbox{Cgrey}{

\begin{Theorem} \label{TLL1}
Let $\mathcal{F}\in C^2(R^{3 \times 3}_{0,{\rm sym}}; \mathbb{R})$ satisfy assumptions $(A1)-(A3)$. 
Then problem (\ref{i1}--\ref{i4}) admits a dissipative solution $[\vc{v}, \mathbb{P}, \mathbb{Q}]$ in $(0,T) \times \tor$ for any initial data
\[
\vc{v}_0 \in L^2(\tor; R^3), \ \Div \vc{v}_0 = 0, \ \mathbb{P}_0 \in L^2(\tor; R^{3 \times 3}_{0, {\rm sym}}),\
\mathbb{Q}_0 \in W^{1,2}(\tor; R^{3 \times 3}_{0, {\rm sym}}).
\]
\end{Theorem}

}



\section{Relative energy and weak-strong uniqueness}
\label{L}

The dissipative solutions introduced in the preceding section may seem rather weak as we have apparently no information about the specific form
of neither the dissipation defect $\mathcal{D}$ nor the correctors $\mathbb{R}^1$, $\mathbb{R}^2$. Nevertheless, we show that a dissipative solution
coincides with the strong solution emanating from the same initial data as long as the latter exists.


\subsection{Relative energy}

We consider the modified energy functional
\[
E(\vv, \mathbb{P}, \mathbb{Q} ) = \intO{ \left[ \frac{1}{2} |\vc{v}|^2 + \frac{1}{2} |\mathbb{P}|^2 + \frac{1}{2} |\Grad \mathbb{Q}|^2 + \mathcal{G}(\mathbb{Q}) \right] }
\]
introduced in (\ref{defi:E}), along with the associated \emph{relative energy functional}
\begin{equation} \label{L2}
\begin{split}
&\mathcal{E} \left( \vc{v}, \mathbb{P}, \mathbb{Q} \ \Big| \tilde {\vc{v}}, \tilde {\mathbb{P}}, \tilde {\mathbb{Q}} \right)\\
&= \frac{1}{2} \intO{ \left[ |\vc{v} - \tv|^2 + |\mathbb{P} - \tP |^2 + | \Grad \mathbb{Q} - \Grad \tQ |^2 \right] }  + \intO{ \mathcal{G} (\mathbb{Q}) -
\partial{\mathcal{G}}(\tQ): (\mathbb{Q} - \tQ) - \mathcal{G}(\tQ)  }\\
&= E(\vc{v}, \mathbb{P}, \mathbb{Q}) + E(\tv, \tP, \tQ ) - \intO{ \left[ \vc{v} \cdot \tv + \mathbb{P}: \tP + \Grad \mathbb{Q} : \Grad \tQ \right] }\\
&- \intO{ \left[ \partial{\mathcal{G}}(\tQ): (\mathbb{Q} - \tQ) + 2 \mathcal{G} (\tQ) \right] }
\end{split}
\end{equation}
defined
for any trio of smooth function $[\tv, \tP, \tQ]$. The functional $\mathcal{E}$ plays a role of a ``distance'' between a solution $[\vc{v}, \mathbb{P},
\mathbb{Q}]$ and a generic triplet $[\tv, \tP, \tQ]$.


\subsection{Relative energy inequality}

Our goal is to derive the relative entropy inequality - an explicit formula for
\[
\left[ \mathcal{E} \left( \vc{v}, \mathbb{P}, \mathbb{Q} \ \Big| \tilde {\vc{v}}, \tilde {\mathbb{P}}, \tilde {\mathbb{Q}} \right) \right]_{t = 0}^{t =
\tau }.
\]
To this aim, we shall directly assume that $[\tv, \tP, \tQ]$ is a strong solution since this simplifies a bit 
the computations. That said, we first use the weak formulation (\ref{WS}) to compute
\begin{equation} \label{L3}
\begin{split}
\left[ \intO{ \vc{v} \cdot \tv } \right]_{t = 0}^{t = \tau} &= \int_0^\tau \intO{
\vv \cdot \partial_t \tv + (\vv \otimes \vv) : \Grad \tv + (\Grad \mathbb{Q} \odot \Grad \mathbb{Q}) : \Grad \tv } \dt  \\
&+ \int_0^\tau \left< \mathbb{R}^1 ; \Grad \tv \right> \dt;
\end{split}
\end{equation}
\begin{equation} \label{L4}
\begin{split}
\left[ \intO{ \mathbb{P} : \tP } \right]_{t = 0}^{t = \tau} &= \int_0^\tau \intO{ \left[
\mathbb{P}: \partial_t \tP  + (\vv \otimes \mathbb{P} ): \Grad \tP  -  \frac{\partial \mathcal{G}(\mathbb{Q}) }{\partial \mathbb{Q} }
: \tP + 2 \Lambda \mathbb{Q}: \tP - \Grad \mathbb{Q} \cdot \Grad \tP \right] } \dt \\
&+ \int_0^\tau \left< \mathbb{R}^2; \Grad \tP \right> \dt;
\end{split}
\end{equation}
\begin{equation} \label{L5}
\begin{split}
\left[ \intO{ \Grad \mathbb{Q} : \Grad \tQ } \right]_{t = 0}^{\tau} &=
- \left[ \intO{ \mathbb{Q} : \Del \tQ  } \right]_{t = 0}^{t = \tau} \\ &= - \int_0^\tau \intO{ \left[
\mathbb{Q} : \partial_t \Del \tQ  + (\vv \otimes \mathbb{Q}) : \Grad \Del \tQ + \mathbb{P}: \Del \tQ \right]} \dt;
\end{split}
\end{equation}
and, finally,
\begin{equation} \label{L6}
\left[ \intO{ \mathbb{Q} : \partial \mathcal{G}(\tQ) } \right]_{t = 0}^{t = \tau} = \int_0^\tau \intO{ \left[
\mathbb{Q} : \partial_t \partial \mathcal{G} (\tQ)  + (\vv \otimes \mathbb{Q}) : \Grad \partial \mathcal{G} (\tQ) + \mathbb{P}: \partial \mathcal{G}
(\tQ) \right]} \dt.
\end{equation}

Using the energy balance (\ref{WSE}), together with the relations (\ref{L3}--\ref{L6})
and the fact that $(\tv,\tP,\tQ)$ is assumed to be a strong solution,
we may compute
\begin{equation} \label{L7}
\begin{split}
&\left[ \mathcal{E} \left( \vc{v}, \mathbb{P}, \mathbb{Q} \ \Big| \tilde {\vc{v}}, \tilde {\mathbb{P}}, \tilde {\mathbb{Q}} \right)
\right]_{t = 0}^{t = \tau} \\
&= \left[ E(\vc{v}, \mathbb{P}, \mathbb{Q}) \right]_{t = 0}^{t = \tau} + \left[ E(\tv, \tP, \tQ )
\right]_{t = 0}^{t = \tau} - \left[ \intO{ \left[ \vc{v} \cdot \tv + \mathbb{P}: \tP + \Grad \mathbb{Q} : \Grad \tQ \right] }
\right]_{t = 0}^{t = \tau} \\
&- \left[ \intO{ \left[ \partial{\mathcal{G}}(\tQ): (\mathbb{Q} - \tQ) + 2 \mathcal{G} (\tQ) \right] } \right]_{t = 0}^{ t = \tau} \\
&= - \mathcal{D} (\tau) + 2\Lambda \int_0^\tau \intO{ \left( \mathbb{Q}: \mathbb{P} + \tQ :\tP \right) } \dt   - \left[ \intO{ \left[ \vc{v} \cdot \tv + \mathbb{P}: \tP + \Grad \mathbb{Q} : \Grad \tQ \right] }
\right]_{t = 0}^{t = \tau} \\
&- \left[ \intO{ \left[ \partial{\mathcal{G}}(\tQ): (\mathbb{Q} - \tQ) + 2 \mathcal{G} (\tQ) \right] } \right]_{t = 0}^{ t = \tau}\\
&= - \mathcal{D} (\tau) + 2 \Lambda \int_0^\tau \intO{ \left[ \mathbb{Q}: \mathbb{P} + \tQ :\tP -  \mathbb{Q} : \tP \right] } \dt \\
&- \int_0^\tau \intO{\left[
\vv \cdot \partial_t \tv + \vv \otimes \vv : \Grad \tv + \Grad \mathbb{Q} \odot \Grad \mathbb{Q}: \Grad \tv \right] } \dt  \\
&- \int_0^\tau \intO{ \left[
\mathbb{P}: \partial_t \tP  + (\vv \otimes \mathbb{P} ): \Grad \tP  -  \frac{\partial \mathcal{G}(\mathbb{Q}) }{\partial \mathbb{Q} }
: \tP  - \Grad \mathbb{Q} \cdot \Grad \tP \right] } \dt \\
&+ \int_0^\tau \intO{ \left[
\mathbb{Q} : \partial_t \Del \tQ  + (\vv \otimes \mathbb{Q}) : \Grad \Del \tQ + \mathbb{P}: \Del \tQ \right]} \dt\\
&-\int_0^\tau \intO{ \left[
\mathbb{Q} : \partial_t \partial \mathcal{G} (\tQ)  + (\vv \otimes \mathbb{Q}) : \Grad \partial \mathcal{G} (\tQ) + \mathbb{P}: \partial \mathcal{G}
(\tQ) \right]} \dt
\\
&- \int_0^\tau \intO{ \partial_t \left( 2 \mathcal{G} (\tQ) - \partial \mathcal{G} (\tQ): \tQ \right) } \dt
- \int_0^\tau \left< \mathbb{R}^2; \Grad \tP \right> \dt - \int_0^\tau \left< \mathbb{R}^1 ; \Grad \tv \right> \dt.
\end{split}
\end{equation}

Using (\ref{DDI}), we may
rewrite (\ref{L7}) in a slightly more concise form obtaining the \emph{relative energy inequality}:
\begin{equation} \label{L8}
\begin{split}
&\left[ \mathcal{E} \left( \vc{v}, \mathbb{P}, \mathbb{Q} \ \Big| \tilde {\vc{v}}, \tilde {\mathbb{P}}, \tilde {\mathbb{Q}} \right)
\right]_{t = 0}^{t = \tau} + \mathcal{D}(\tau) \\
&\leq 2\Lambda \int_0^\tau \intO{ \left[ \mathbb{Q}: \mathbb{P} + \tQ :\tP -  \mathbb{Q} : \tP \right] } \dt \\
&- \int_0^\tau \intO{\left[
\vv \cdot \partial_t \tv + \vv \otimes \vv : \Grad \tv + \Grad \mathbb{Q} \odot \Grad \mathbb{Q}: \Grad \tv \right] } \dt  \\
&- \int_0^\tau \intO{ \left[
\mathbb{P}: \partial_t \tP  + (\vv \otimes \mathbb{P} ): \Grad \tP  -  \frac{\partial \mathcal{G}(\mathbb{Q}) }{\partial \mathbb{Q} }
: \tP  - \Grad \mathbb{Q} \cdot \Grad \tP \right] } \dt \\
&+ \int_0^\tau \intO{ \left[
\mathbb{Q} : \partial_t \Del \tQ  + (\vv \otimes \mathbb{Q}) : \Grad \Del \tQ + \mathbb{P}: \Del \tQ \right]} \dt\\
&-\int_0^\tau \intO{ \left[
\mathbb{Q} : \partial_t \partial \mathcal{G} (\tQ)  + (\vv \otimes \mathbb{Q}) : \Grad \partial \mathcal{G} (\tQ) + \mathbb{P}: \partial \mathcal{G}
(\tQ) \right]} \dt
\\
&- \int_0^\tau \intO{ \partial_t \left( 2 \mathcal{G} (\tQ) - \partial \mathcal{G} (\tQ): \tQ \right) } \dt
+ c \int_0^\tau \left( \| \Grad \tP \|_{C(\tor)} + \|\Grad \tv \|_{C(\tor)} \right) \mathcal{D}(\cdot) \dt,
\end{split}
\end{equation}
holding for any strong solution $[\tv, \tP, \tQ]$.


\subsection{Weak-strong uniqueness}

Our ultimate goal is to show that any dissipative solution necessarily coincides with a strong solution originating from the same initial data on
the existence interval of the latter. A simple idea is to take the strong solution $[\tv, \tP, \tQ]$ as ``test functions'' in the relative energy inequality
(\ref{L8}) and to use a Gronwall-type argument. This will be done in several steps.


\subsubsection{Step 1 - velocity}

Using the fact that the velocity $\tv$ satisfies (\ref{i1}), (\ref{i2}) we deduce
\[
\begin{split}
- \int_0^\tau &\intO{\left[
\vv \cdot \partial_t \tv + \vv \otimes \vv : \Grad \tv + \Grad \mathbb{Q} \odot \Grad \mathbb{Q}: \Grad \tv \right] } \dt\\
& = \int_0^\tau \intO{\left[
\vv \cdot \left(\tv \cdot \Grad\right) \tv - \vv \otimes \vv :\Grad \tv +
\Div \left( \Grad \tQ \odot \Grad \tQ \right) \cdot \vc{v}
- \Grad \mathbb{Q} \odot \Grad \mathbb{Q}: \Grad \tv \right] } \dt\\
& = -\int_0^\tau \intO{ (\vv - \tv) \cdot \Grad \tv \cdot (\vv - \tv) } \dt \\
&+ \int_0^\tau \intO{ \left[ \Div \left( \Grad \tQ \odot \Grad \tQ \right) \cdot \vc{v}
- \Grad \mathbb{Q} \odot \Grad \mathbb{Q}: \Grad \tv \right] } \dt.
\end{split}
\]
Consequently, inequality (\ref{L8}) reads
\begin{equation} \label{L9}
\begin{split}
&\left[ \mathcal{E} \left( \vc{v}, \mathbb{P}, \mathbb{Q} \ \Big| \tilde {\vc{v}}, \tilde {\mathbb{P}}, \tilde {\mathbb{Q}} \right)
\right]_{t = 0}^{t = \tau} + \mathcal{D}(\tau) \\
&\leq 2\Lambda \int_0^\tau \intO{ \left[ \mathbb{Q}: \mathbb{P} + \tQ :\tP -  \mathbb{Q} : \tP \right] } \dt
- \int_0^\tau \intO{ (\vv - \tv) \cdot \Grad \tv \cdot (\vv - \tv) } \dt \\
&+ \int_0^\tau \intO{ \left[ \Div \left( \Grad \tQ \odot \Grad \tQ \right) \cdot \vc{v}
- \Grad \mathbb{Q} \odot \Grad \mathbb{Q}: \Grad \tv \right] } \dt\\
&- \int_0^\tau \intO{ \left[
\mathbb{P}: \partial_t \tP  + (\vv \otimes \mathbb{P} ): \Grad \tP  -  \frac{\partial \mathcal{G}(\mathbb{Q}) }{\partial \mathbb{Q} }
: \tP  - \Grad \mathbb{Q} \cdot \Grad \tP \right] } \dt \\
&+ \int_0^\tau \intO{ \left[
\mathbb{Q} : \partial_t \Del \tQ  + (\vv \otimes \mathbb{Q}) : \Grad \Del \tQ + \mathbb{P}: \Del \tQ \right]} \dt\\
&-\int_0^\tau \intO{ \left[
\mathbb{Q} : \partial_t \partial \mathcal{G} (\tQ)  + (\vv \otimes \mathbb{Q}) : \Grad \partial \mathcal{G} (\tQ) + \mathbb{P}: \partial \mathcal{G}
(\tQ) \right]} \dt
\\
&- \int_0^\tau \intO{ \partial_t \left( 2 \mathcal{G} (\tQ) - \partial \mathcal{G} (\tQ): \tQ \right) } \dt
+ c \int_0^\tau \left( \| \Grad \tP \|_{C(\tor)} + \|\Grad \tv \|_{C(\tor)} \right) \mathcal{D}(\cdot) \dt.
\end{split}
\end{equation}
%


\subsubsection{Step 2 - $\mathbb{P}$ tensor}

As $\tP$ satisfies (\ref{i4}), we get
\[
\begin{split}
- \int_0^\tau &\intO{ \left[
\mathbb{P}: \partial_t \tP  + (\vv \otimes \mathbb{P} ): \Grad \tP  -  \frac{\partial \mathcal{G}(\mathbb{Q}) }{\partial \mathbb{Q} }
: \tP  - \Grad \mathbb{Q} \cdot \Grad \tP \right] } \dt \\
&= \int_0^\tau \intO{ \left[
\mathbb{P}:(\tv \cdot \Grad) \tP + \mathbb{P} : \partial \mathcal{F}(\tQ) - \Del \tQ : \mathbb{P}  - (\vv \otimes \mathbb{P} ): \Grad \tP \right] } \dt \\
&+  \int_0^\tau \intO{ \left[ \frac{\partial \mathcal{G}(\mathbb{Q}) }{\partial \mathbb{Q} }
: \tP  + \Grad \mathbb{Q} \cdot \Grad \tP \right] } \dt\\
&=
\int_0^\tau \intO{ \left[ \left( (\tv - \vv) \cdot \Grad \tP\right) : (\mathbb{P} - \tP )
 + \mathbb{P} : \partial \mathcal{G}(\tQ) - 2\Lambda \tQ : \mathbb{P} - \Del \tQ : \mathbb{P}  \right] } \dt \\
&+  \int_0^\tau \intO{ \left[ {\partial \mathcal{G}(\mathbb{Q}) }: \tP  + \Grad \mathbb{Q} \cdot \Grad \tP \right] } \dt.
\end{split}
\]
Thus we may rewrite (\ref{L9}) as
\begin{equation} \label{L10}
\begin{split}
&\left[ \mathcal{E} \left( \vc{v}, \mathbb{P}, \mathbb{Q} \ \Big| \tilde {\vc{v}}, \tilde {\mathbb{P}}, \tilde {\mathbb{Q}} \right)
\right]_{t = 0}^{t = \tau} + \mathcal{D}(\tau) \\
&\leq 2\Lambda \int_0^\tau \intO{ (\mathbb{Q} - \tQ): (\mathbb{P} - \tP) } \dt
- \int_0^\tau \intO{ (\vv - \tv) \cdot \Grad \tv \cdot (\vv - \tv) } \dt \\
&+
\int_0^\tau \intO{  (\tv - \vv) \cdot \Grad \tP \cdot (\mathbb{P} - \tP ) } \dt \\
&+  \int_0^\tau \intO{ \left[ {\partial \mathcal{G}(\mathbb{Q}) }
: \tP  + \Grad \mathbb{Q} \cdot \Grad \tP \right] } \dt\\
&+ \int_0^\tau \intO{ \left[ \Div \left( \Grad \tQ \odot \Grad \tQ \right) \cdot \vc{v}
- \Grad \mathbb{Q} \odot \Grad \mathbb{Q}: \Grad \tv \right] } \dt\\
&+ \int_0^\tau \intO{ \left[
\mathbb{Q} : \partial_t \Del \tQ  + (\vv \otimes \mathbb{Q}) : \Grad \Del \tQ  \right]} \dt\\
&-\int_0^\tau \intO{ \left[
\mathbb{Q} : \partial_t \partial \mathcal{G} (\tQ)  + (\vv \otimes \mathbb{Q}) : \Grad \partial \mathcal{G} (\tQ)  \right]} \dt
\\
&- \int_0^\tau \intO{ \partial_t \left( 2 \mathcal{G} (\tQ) - \partial \mathcal{G} (\tQ): \tQ \right) } \dt
+ c \int_0^\tau \left( \| \Grad \tP \|_{C(\tor)} + \|\Grad \tv \|_{C(\tor)} \right) \mathcal{D}(\cdot) \dt.
\end{split}
\end{equation}


\subsubsection{Step 3 - $\mathbb{Q}$ tensor}

First of all, let us compute
\[
\begin{split}
\int_0^\tau &\intO{ \left[
\mathbb{Q} : \partial_t \Del \tQ  + (\vv \otimes \mathbb{Q}) : \Grad \Del \tQ  \right]} \dt \\
&=\int_0^\tau \intO{ \left[-
\Grad \mathbb{Q} \cdot \partial_t \Grad \tQ  + (\vv \otimes \mathbb{Q}) : \Grad \Del \tQ  \right]} \dt\\
&=\int_0^\tau \intO{ \left[\Grad \Div (\tv \otimes \tQ ) \cdot \Grad \mathbb{Q} - \Grad \tP \cdot \Grad \mathbb{Q}
+ (\vv \otimes \mathbb{Q}) : \Grad \Del \tQ  \right]} \dt.
\end{split}
\]
Here, in deducing the latter equality, we used (\ref{i3}), written for the strong solution 
$\tQ$, and applied to it the operator $\Grad$. Similarly, 
multiplying (\ref{i3}) in the tensor sense by 
(the fourth order tensor) $-\frac{\partial^2 \mathcal{G}}{\partial^2 \mathbb{Q}}$, we infer
\[
\begin{split}
-\int_0^\tau &\intO{ \left[
\mathbb{Q} : \partial_t \partial \mathcal{G} (\tQ)  + (\vv \otimes \mathbb{Q}) : \Grad \partial \mathcal{G} (\tQ)  \right]} \dt \\
= \int_0^\tau &\intO{ \left[
\mathbb{Q} : \left( \tv \cdot \Grad \partial \mathcal{G} (\tQ) \right) 
 - \tP : \left( \partial^2 \mathcal{G}(\tQ) \mathbb{Q} \right) -  (\vv \otimes \mathbb{Q}) :\Grad \partial \mathcal{G} (\tQ)  \right]} \dt.
\end{split}
\]
Then, noting also that 
\[
\begin{split}
 \int_0^\tau &\intO{ \left[
  \mathbb{Q} : \big( \tv \cdot \Grad \partial \mathcal{G} (\tQ) \big) 
  -  (\vv \otimes \mathbb{Q}) :\Grad \partial \mathcal{G} (\tQ)  \right] } \dt\\
 = \int_0^\tau &\intO{ \left[
   \mathbb{Q} : \big( (\tv - \vc{v}) \cdot \Grad \partial \mathcal{G} (\tQ) \big) \right] } \dt \\
 = \int_0^\tau &\intO{ \left[
   ( \mathbb{Q} - \tQ ) : \big( (\tv - \vc{v}) \cdot \Grad \partial \mathcal{G} (\tQ) \big) 
   + \tQ : \big( (\tv - \vc{v}) \cdot \Grad \partial \mathcal{G} (\tQ) \big) \right] } \dt,
\end{split}
\]
we see that the relative energy inequality (\ref{L10}) 
may be rewritten in the form
\begin{equation} \label{L11}
\begin{split}
&\left[ \mathcal{E} \left( \vc{v}, \mathbb{P}, \mathbb{Q} \ \Big| \tilde {\vc{v}}, \tilde {\mathbb{P}}, \tilde {\mathbb{Q}} \right)
\right]_{t = 0}^{t = \tau} + \mathcal{D}(\tau) \\
&\leq 2\Lambda \int_0^\tau \intO{ (\mathbb{Q} - \tQ): (\mathbb{P} - \tP) } \dt
- \int_0^\tau \intO{ (\vv - \tv) \cdot \Grad \tv \cdot (\vv - \tv) } \dt \\
&+
\int_0^\tau \intO{  (\tv - \vv) \cdot \Grad \tP \cdot (\mathbb{P} - \tP ) } \dt +
\int_0^\tau \intO{ (\tv - \vv) \cdot \Grad \partial  \mathcal{G}(\tQ)  \cdot (\mathbb{Q} - \tQ) } \dt \\
&+  \int_0^\tau \intO{  \tP : \left( \partial \mathcal{G}(\mathbb{Q})  - \partial^2 \mathcal{G}(\tQ) ( \mathbb{Q} -
\tQ) - \partial \mathcal{G} (\tQ)  \right) } \dt\\
&+ \int_0^\tau \intO{ \left[ \Div \left( \Grad \tQ \odot \Grad \tQ \right) \cdot \vc{v}
- \Grad \mathbb{Q} \odot \Grad \mathbb{Q}: \Grad \tv \right] } \dt\\
&+
\int_0^\tau \intO{ \left[\Grad \Div (\tv \otimes \tQ ) \cdot \Grad \mathbb{Q}
+ (\vv \otimes \mathbb{Q}) : \Grad \Del \tQ  \right]} \dt\\
& + \int_0^\tau \intO{ \left[ \tP : \partial \mathcal{G}(\tQ) - \tP : \partial^2 \mathcal{G}(\tQ): \tQ \right] } \dt
- \int_0^\tau \intO{ \partial_t \left( 2 \mathcal{G} (\tQ) - \partial \mathcal{G} (\tQ): \tQ \right) } \dt \\
&+ c \int_0^\tau \left( \| \Grad \tP \|_{C(\tor)} + \|\Grad \tv \|_{C(\tor)} \right) \mathcal{D}(\cdot) \dt.
\end{split}
\end{equation}


\subsubsection{Conclusion}

As $\tQ$ satisfies the transport equation (\ref{i3}), we easily deduce that
\[
\int_0^\tau \intO{ \left[ \tP : \partial \mathcal{G}(\tQ) - \tP : \partial^2 \mathcal{G}(\tQ): \tQ \right] } \dt
- \int_0^\tau \intO{ \partial_t \left( 2 \mathcal{G} (\tQ) - \partial \mathcal{G} (\tQ): \tQ \right) } \dt = 0.
\]
Moreover, after a straightforward manipulation,
\[
\begin{split}
\int_0^\tau &\intO{ \left[ \Div \left( \Grad \tQ \odot \Grad \tQ \right) \cdot \vc{v}
- \Grad \mathbb{Q} \odot \Grad \mathbb{Q}: \Grad \tv \right] } \dt\\
&+
\int_0^\tau \intO{ \left[\Grad \Div (\tv \otimes \tQ ) \cdot \Grad \mathbb{Q}
+ (\vv \otimes \mathbb{Q}) : \Grad \Del \tQ  \right]} \dt\\
&=\int_0^\tau \intO{ \left[ \Div \left( \Grad \tQ \odot \Grad \tQ \right) \cdot (\vc{v} - \tv )
- \Grad \mathbb{Q} \odot \Grad \mathbb{Q}: \Grad \tv \right] } \dt\\
&+
\int_0^\tau \intO{ \left[\Grad \Div (\tv \otimes \tQ ) \cdot ( \Grad \mathbb{Q} - \Grad \tQ)
- \vv \cdot \Grad \mathbb{Q} \cdot \Del \tQ  \right]} \dt\\
&=\int_0^\tau \intO{ \left[ \Grad \tQ \cdot \Del \tQ \cdot (\vc{v} - \tv )
- \Grad \mathbb{Q} \odot \Grad \mathbb{Q}: \Grad \tv \right] } \dt\\
&+
\int_0^\tau \intO{ \left[\Grad \Div (\tv \otimes \tQ ) \cdot ( \Grad \mathbb{Q} - \Grad \tQ)
- \vv \cdot \Grad \mathbb{Q} \cdot \Del \tQ  \right]} \dt\\
&=\int_0^\tau \intO{ \left[ (\Grad \tQ - \Grad \mathbb{Q})  \cdot \Del \tQ \cdot (\vc{v} - \tv )
- \Grad \mathbb{Q} \odot \Grad \mathbb{Q}: \Grad \tv \right] } \dt\\
&+
\int_0^\tau \intO{ \left[\Grad (\tv \cdot \Grad \tQ ) \cdot ( \Grad \mathbb{Q} - \Grad \tQ)
- \tv \cdot \Grad \mathbb{Q} \cdot \Del \tQ  \right]} \dt\\
&=\int_0^\tau \intO{ \left[ (\Grad \tQ - \Grad \mathbb{Q})  \cdot \Del \tQ \cdot (\vc{v} - \tv ) -
(\Grad \mathbb{Q} - \Grad \tQ) \cdot \Grad \tv \cdot (\Grad \mathbb{Q} - \Grad \tQ) \right] } \dt \\
&
- \int_0^\tau \intO{ \Grad \tQ \odot \Grad \mathbb{Q} : \Grad \tv  } \dt
+
\int_0^\tau \intO{ \left[ \tv \cdot \Grad \Grad \tQ  \cdot  \Grad \mathbb{Q}
- \tv \cdot \Grad \mathbb{Q} \cdot \Del \tQ  \right]} \dt\\
&=\int_0^\tau \intO{ \left[ (\Grad \tQ - \Grad \mathbb{Q})  \cdot \Del \tQ \cdot (\vc{v} - \tv ) -
(\Grad \mathbb{Q} - \Grad \tQ) \cdot \Grad \tv \cdot (\Grad \mathbb{Q} - \Grad \tQ) \right] } \dt,
\end{split}
\]
where we have used solenoidality of $\tv$ to observe that
\[
\intO{ \left[ - \Grad \tQ \odot \Grad \mathbb{Q} : \Grad \tv
+ \tv \cdot \Grad \Grad \tQ  \cdot  \Grad \mathbb{Q}
- \tv \cdot \Grad \mathbb{Q} \cdot \Del \tQ  \right]}  = 0
\]
Accordingly, the relative energy inequality takes the final form
\begin{equation} \label{L12}
\begin{split}
& \mathcal{E} \left( \vc{v}, \mathbb{P}, \mathbb{Q} \ \Big| \tilde {\vc{v}}, \tilde {\mathbb{P}}, \tilde {\mathbb{Q}} \right)
(\tau) + \mathcal{D}(\tau) \\
&\leq 2\Lambda \int_0^\tau \intO{ (\mathbb{Q} - \tQ): (\mathbb{P} - \tP) } \dt
-\int_0^\tau \intO{ (\vv - \tv) \cdot \Grad \tv \cdot (\vv - \tv) } \dt \\
&+
\int_0^\tau \intO{  (\tv - \vv) \cdot \Grad \tP \cdot (\mathbb{P} - \tP ) } \dt +
\int_0^\tau \intO{ (\tv - \vv) \cdot \Grad \partial  \mathcal{G}(\tQ)  \cdot (\mathbb{Q} - \tQ) } \dt \\
&+  \int_0^\tau \intO{  \tP : \left( \partial \mathcal{G}(\mathbb{Q})  - \partial^2 \mathcal{G}(\tQ) ( \mathbb{Q} -
\tQ) - \partial \mathcal{G} (\tQ)  \right) } \dt\\
&+\int_0^\tau \intO{ \left[ (\Grad \tQ - \Grad \mathbb{Q})  \cdot \Del \tQ \cdot (\vc{v} - \tv ) -
(\Grad \mathbb{Q} - \Grad \tQ) \cdot \Grad \tv \cdot (\Grad \mathbb{Q} - \Grad \tQ) \right] } \dt\\
&+ c \int_0^\tau \left( \| \Grad \tP \|_{C(\tor)} + \|\Grad \tv \|_{C(\tor)} \right) \mathcal{D}(\cdot) \dt.
\end{split}
\end{equation}

Applying Gronwall's lemma we get the desired conclusion:

\Cbox{Cgrey}{

\begin{Theorem} \label{TWS}

Under the hypotheses of Theorem \ref{TLL1},
let the initial data enjoy the regularity properties (\ref{i5a}). Let $[\vv, \mathbb{P}, \mathbb{Q}]$ be a dissipative solution of
problem (\ref{i1}--\ref{i4}), (\ref{i5a}) and let $[\tv, \tP, \tQ]$ a the strong solution of the same problem
belonging to the regularity class specified in Theorem \ref{TE1}
in the space-time cylinder $(0,T) \times \tor$.

Then
\[
\vv = \tv,\ \mathbb{P} = \tP,\ \mathbb{Q} = \tQ \ \mbox{a.a. in}\ (0,T) \times \tor.
\]

\end{Theorem}

}

\begin{Remark} \label{KK1}
The assumptions concerning regularity of the strong solution are not optimal and may be 
relaxed. 
\end{Remark}
Combining Theorem \ref{TWS} with the local existence result established in Theorem \ref{TE1} we immediately get the following 
corollary:
\begin{Corollary} \label{cor}
Let $[\vv, \mathbb{P}, \mathbb{Q}]$ be a dissipative solution of problem (\ref{i1}--\ref{i4}) in $(0,T) \times \tor$ enjoying the regularity
specified in\/ {\rm Theorem \ref{TE1}}.

Then $[\vv, \mathbb{P}, \mathbb{Q}]$ is a strong solution, in particular, the dissipation defect $\mathcal{D}$ as well as the
defect measures $\mathbb{R}^1$, $\mathbb{R}^2$ vanish identically in $[0,T] \times \tor$.
\end{Corollary}

\def\cprime{$'$} \def\ocirc#1{\ifmmode\setbox0=\hbox{$#1$}\dimen0=\ht0
  \advance\dimen0 by1pt\rlap{\hbox to\wd0{\hss\raise\dimen0
  \hbox{\hskip.2em$\scriptscriptstyle\circ$}\hss}}#1\else {\accent"17 #1}\fi}


\section*{Acknowledgement}

The research of E.F.~leading to these results has received funding from the
European Research Council under the European Union's Seventh
Framework Programme (FP7/2007-2013)/ ERC Grant Agreement
320078. The Institute of Mathematics of the Academy of Sciences of
the Czech Republic is supported by RVO:67985840.
A.Z.~was partially supported by a grant of the Romanian National Authority for 
Scientific Research and Innovation, CNCS-UEFISCDI, project number PN-II-RU-TE-2014-4-0657.
The financial support of the FP7-IDEAS-ERC-StG \#256872
(EntroPhase) is gratefully acknowledged by the authors.
E.R.~and G.S.~have been partially supported by GNAMPA, the 
``Gruppo Nazionale per l'Analisi Matematica, la Probabilit\`a e le loro Applicazioni'' 
of INdAM  and the IMATI -- C.N.R. Pavia.



\end{document}